\documentclass{article}

\usepackage{graphicx}
\usepackage{lmodern}
\usepackage{amsmath,amssymb,amsthm,url}

\newtheorem{theorem}{Theorem}[section]
\newtheorem{proposition}[theorem]{Proposition}
\newtheorem{lemma}[theorem]{Lemma}

\newcommand{\R}{\mathbb R}
\newcommand{\conv}{\operatorname{conv}}
\newcommand{\spanof}{\operatorname{span}}

\title{Antipodal paths: a short proof by counting intersections}

\author{Valjakas Djaljapayan\thanks{%
Department of Mathematics, National Central University,
No.~300, Zhongda Road, Zhongli District,
Taoyuan City 320317, Taiwan.
Email: \texttt{115281004@cc.ncu.edu.tw}.
}}

\date{}

\begin{document}

\maketitle

\begin{abstract}
We give a short elementary proof of Norine's antipodal-coloring
theorem: every red--blue edge coloring of a hypercube in which
antipodal edges have opposite colors contains a monochromatic
path joining antipodal vertices.
\end{abstract}

\section{Introduction}

For an integer $d\ge1$, the cube graph $Q_d$ has vertex set
$\{0,1\}^d$; two vertices are adjacent when they differ in exactly
one coordinate. For $x=(x_1,\ldots,x_d)\in\{0,1\}^d$, write
\[
 \bar x=(1-x_1,\ldots,1-x_d).
\]
The vertices $x$ and $\bar x$ are called antipodal, and the
antipodal edge of $xy$ is $\bar x\bar y$.
A red--blue edge coloring is called antipodal if every edge has
the opposite color to its antipodal edge.

Norine's conjecture, posed in 2008 and discussed by Feder and
Subi~\cite{FS}, asks whether every antipodal edge coloring of a
cube contains a monochromatic path joining antipodal vertices.
Wu and Yang~\cite{WY} proved the following statement.

\begin{theorem}[Norine's antipodal-coloring theorem]
\label{thm:main}
Let $n\ge2$. Suppose that the edges of $Q_n$ are colored red or blue,
and every edge has the opposite color to its antipodal edge.
Then some vertex $y$ is joined to $\bar y$ by a red path.
\end{theorem}

Taking the antipodal image of a red path gives a blue path between
the same two vertices, so either color can be specified in the
conclusion.
A stronger question asks for a monochromatic antipodal
\emph{geodesic}: a path that changes every coordinate exactly once
and therefore has length $n$.
Theorem~\ref{thm:main} concerns the existence of an antipodal path
of unrestricted length.

Several earlier results addressed structural restrictions on the
coloring. Feder and Subi~\cite{FS} proved that an arbitrary
red--blue edge coloring has a monochromatic antipodal path if
there is no square whose edge colors alternate red, blue, red,
blue. Solt\'esz~\cite{Soltesz} weakened this restriction and studied
related questions for graphs equipped with an automorphism.
Feder and Subi also formulated the stronger conjecture that
every red--blue edge coloring of $Q_n$ has an antipodal path with
at most one color change. For an antipodal coloring, such a path
can be converted into a monochromatic antipodal walk by taking
the antipodal image of one of its two monochromatic parts.
Removing repeated vertices then gives the required path.

Finite-dimensional cases were established by a combination of
structural arguments and computation.
West and Wise~\cite{WW} proved the stronger geodesic conclusion
for $2\le n\le6$.
Frankston and Scheinerman~\cite{FrankstonScheinerman} used
Boolean satisfiability (SAT) solvers to verify the conjecture
for $n=7$.
Kirchweger, Peitl, Subercaseaux, and Szeider~\cite{KPSS}
subsequently extended the verification to $n=8$, including
the geodesic version.

Another line of work studies long monochromatic geodesics and
antipodal geodesics with few color changes.
Leader and Long~\cite{LL} proved that a subgraph of $Q_n$ with
average degree $a$ contains a geodesic of length at least $a$.
In particular, every red--blue edge coloring contains a
monochromatic geodesic of length at least $\lceil n/2\rceil$.
For antipodal geodesics in arbitrary red--blue edge colorings,
Dvo\v{r}\'ak~\cite{Dvorak} obtained an upper bound of
$(3/8+o(1))n$ on the number of color changes.
Kirchweger et al.~\cite{KPSS} improved this to $5n/16+6$,
and Hollom~\cite{Hollom} subsequently proved the sublinear bound
\[
 \left(\frac{\pi}{2}+o(1)\right)\sqrt n.
\]

Topological methods provide a further approach.
D\v{z}avoronok~\cite{Dzavoronok} used a criterion for triangulated
$2$-complexes to establish monochromatic antipodal connectivity
for certain classes of cube colorings.
The proof of Wu and Yang~\cite{WY} treats all antipodal colorings
by deriving an equivariant chain map from an assumed
counterexample and ruling it out with a chain-level
Borsuk--Ulam obstruction.

This paper is motivated by the search for a short proof in which
the parity argument is expressed through explicit finite counts.
We retain the component labeling and parity framework of
Wu and Yang~\cite{WY}, and formulate the argument using convex
hulls of path labels and nested linear subspaces.
Under an assumed counterexample, the first intersection count
is one and the last is zero.
We show that consecutive counts have the same parity.
The local identity reduces to the fact that a nondegenerate
closed line segment has two endpoints, whereas a ray has one.
Summing these identities, pairing coordinate orders, and counting
face incidences gives the contradiction.
All of these steps are proved directly using elementary linear
algebra, convexity, and counting.

The paper is organised as follows.
Section~\ref{Com} introduces the component labels and intersection
counts, proves the rank bound needed to handle convex hulls
containing the origin, and establishes the local parity identity
by counting the endpoints of a line segment or a ray.
Section~\ref{local} combines this identity with cancellations among
monotone paths and a count of face incidences to show that
consecutive intersection counts have the same parity,
completing the proof of Theorem~\ref{thm:main}.

\section{Component labels and the contradiction assumption}
\label{Com}

Assume, for a contradiction, that no vertex is joined to its
antipode by a red path.
Let $C_1,\ldots,C_r$ be the vertex sets of the connected components
of the spanning subgraph of $Q_n$ consisting of the red edges,
including isolated vertices. Define
\[
 p:\{0,1\}^n\longrightarrow\{1,\ldots,r\},
 \qquad
 p(y)=a\quad\text{if }y\in C_a.
\]
This map is well defined because every vertex belongs to exactly
one of the sets $C_1,\ldots,C_r$.
Thus $p(y)$ is the index of the red component containing $y$.
In particular, $p(y)=p(z)$ if and only if $y,z$ are joined by a
red path, allowing a path of length zero when $y=z$.
The assumption says exactly that
\begin{equation}\label{eq:distinct}
 p(y)\ne p(\bar y)\qquad(y\in\{0,1\}^n).
\end{equation}

Choose $K\ge\max\{n,r,2\}$ and put $[K]=\{1,\ldots,K\}$.
Let $\pi:\{0,1\}^K\to\{0,1\}^n$ be projection onto the first
$n$ coordinates.
The additional coordinates allow both coordinate indices and
component numbers to belong to $[K]$.
Label a vertex $x$ of $Q_K$ by the ordered pair
\[
 q(x)=\bigl(p(\pi x),p(\overline{\pi x})\bigr).
\]
These labels satisfy the following three facts:
\begin{enumerate}
 \item The two entries of $q(x)$ are distinct, by
 \eqref{eq:distinct}.
 \item If $q(x)=(a,b)$, then $q(\bar x)=(b,a)$.
 \item Adjacent vertices have the same first entry or the same
 second entry.
\end{enumerate}
Indeed, an original red edge preserves the first entry of $q$.
An original blue edge preserves the second entry, because its
antipodal edge is red.
An edge in an additional coordinate preserves both entries.

Let $e_1,\ldots,e_K$ be the standard basis of $\R^K$, and set
\[
 W=\left\{z\in\R^K:\sum_{a=1}^Kz_a=0\right\},
 \qquad d=K-1.
\]
If $q(x)=(a,b)$, define $v(x)=e_a-e_b\in W$. Thus
\begin{equation}\label{eq:opposite}
 v(x)\ne0,\qquad v(\bar x)=-v(x).
\end{equation}
All these vectors belong to the finite set
\[
 \mathcal V=\{e_a-e_b:a,b\in[K],\ a\ne b\}.
\]
We now construct counts of intersections with nested subspaces
\[
 W=H_0\supset H_1\supset\cdots\supset H_d=\{0\}.
\]
Only intersections away from $0$ will be counted.
This makes the beginning and end of the argument transparent:
the initial count is one, whereas the final count is zero.
The proof consists of showing that the parity cannot change
between them.

\subsection{What is counted, and why the endpoint counts differ}

For a finite set of vectors, its convex hull is the set of their
linear combinations with nonnegative coefficients summing to one.
Its linear span is the set of all their linear combinations.

Choose linear maps $\ell_1,\ldots,\ell_d:W\to\R$ such that,
for every independent tuple $u_1,\ldots,u_h$ from
$\mathcal V$, with $1\le h\le d$,
\begin{equation}\label{eq:forms}
 \det\bigl(\ell_a(u_b)\bigr)_{1\le a,b\le h}\ne0.
\end{equation}
Here and below, independent means linearly independent.
Such maps exist.
Fix a basis of $W$ and express the maps by their coefficients
in that basis.
There are finitely many determinant conditions.
Each determinant is a polynomial that is not identically zero:
for a given independent tuple, one can choose linear maps
whose values on that tuple form the identity matrix.
The product of these polynomials is therefore not identically
zero, so one can choose coefficients at which every determinant
is nonzero.

We record the precise property supplied by this choice.
If $U$ is an $r$-dimensional subspace spanned by vectors from
$\mathcal V$, then
\begin{equation}\label{eq:projection}
 u\longmapsto\bigl(\ell_1(u),\ldots,\ell_h(u)\bigr)
 \quad\text{has rank }\min\{h,r\}\text{ on }U
 \qquad(0\le h\le d).
\end{equation}
To check this, choose a basis of $U$ from $\mathcal V$ and apply
\eqref{eq:forms} to the appropriate number of basis vectors.
For $h=0$ the map is zero.
In particular, it is injective on $U$ when $h\ge r$.

Define
\[
 H_i=\{u\in W:\ell_1(u)=\cdots=\ell_i(u)=0\},
 \qquad 0\le i\le d,
\]
with $H_0=W$.
These subspaces have dimensions $d-i$, so $H_d=\{0\}$.

An $i$-face of $Q_K$ is obtained by fixing $K-i$ coordinates and
allowing the remaining $i$ coordinates to vary.
A monotone path in this face starts with all free coordinates
equal to $0$ and changes each free coordinate once from $0$
to $1$.
Let $\mathcal P(F)$ be the set of such paths in a face $F$.
It has $i!$ elements.
For a $0$-face there is one path, consisting of its single vertex.

A path $P=(x_0,\ldots,x_i)$ determines the convex hull
\[
 S(P)=\conv\{v(x_0),\ldots,v(x_i)\}\subseteq W.
\]
Let $X_i$ be the family of $i$-faces whose fixed bits, read in
increasing coordinate order, are
\[
 0,1,0,1,\ldots.
\]
Choosing the free coordinates determines exactly one such face,
so $|X_i|=\binom Ki$.
The central count is
\begin{equation}\label{eq:counts}
 c_i=\#\left\{(F,P):
 F\in X_i,\ P\in\mathcal P(F),\
 S(P)\cap\bigl(H_i\setminus\{0\}\bigr)\ne\varnothing
 \right\}.
\end{equation}
The symbol $\#$ denotes the number of elements of a finite set.
In other words, $c_i$ is the number of pairs $(F,P)$ with
$F\in X_i$ and $P\in\mathcal P(F)$ such that the convex hull
of the vectors $v(x)$, as $x$ ranges over the vertices of $P$,
contains a nonzero point of $H_i$ for $0\le i\le d$.

The endpoint values now follow directly from the definition.
The family $X_0$ consists of the single vertex
$x_*=(0,1,0,1,\ldots)$.
Since $H_0=W$ and $v(x_*)\ne0$, we have $c_0=1$.
Since $H_d=\{0\}$, there is no nonzero intersection to count
at the other end, then $c_d=0$.
It therefore remains to prove just one assertion:
\begin{equation}\label{eq:goal}
 c_{i-1}\equiv c_i\pmod2\qquad(1\le i\le d).
\end{equation}
If we establish \eqref{eq:goal}, then
\[
 1=c_0\equiv c_1\equiv\cdots\equiv c_d=0\pmod2,
\]
a contradiction.
The remainder of the proof establishes \eqref{eq:goal}.

\subsection{What the edge condition implies about the origin}

We first isolate the consequence of the condition on adjacent
labels.
Its purpose is to control the cases in which a path's
convex hull contains the origin.

\begin{lemma}\label{lem:rank}
Let $(a_1,b_1),\ldots,(a_m,b_m)$ be pairs in $[K]^2$, with
$a_h\ne b_h$.
Suppose consecutive pairs have the same first entry
or the same second entry, and put $u_h=e_{a_h}-e_{b_h}$.
If $0\in\conv\{u_1,\ldots,u_m\}$, then
\[
 \dim\spanof\{u_1,\ldots,u_m\}\le m-2.
\]
\end{lemma}

\begin{proof}
Let $A=\{a_1,\ldots,a_m\}$ and $B=\{b_1,\ldots,b_m\}$.
Reading the pairs in order, the first pair contributes two
to $|A|+|B|$, and each subsequent pair contributes at most one,
because it shares its first or second entry with the preceding
pair.
Thus $|A|+|B|\le m+1$.

Choose $\lambda_1,\ldots,\lambda_m\ge0$ such that
$\sum_h\lambda_h=1$ and $\sum_h\lambda_hu_h=0$.
For any $h$ with $\lambda_h>0$, comparing the $a_h$ and $b_h$
coordinates gives $a_h\in B$ and $b_h\in A$.
Since $a_h\ne b_h$, we obtain $|A\cap B|\ge2$, and therefore
\[
 |A\cup B|=|A|+|B|-|A\cap B|\le m-1.
\]

All $u_h$ have zero coordinates outside $A\cup B$
and have coordinate sum zero, so
\[
 \dim\spanof\{u_1,\ldots,u_m\}
 \le |A\cup B|-1\le m-2.
\]
\end{proof}

For a $j$-step monotone path, this says
\begin{equation}\label{eq:origin}
 0\in S(P)\quad\Longrightarrow\quad \dim S(P)\le j-1.
\end{equation}
Indeed, if a convex hull contains $0$, its affine span is the
linear span of its vertices.
Here $\dim S(P)$ means the dimension of that affine span.
Thus a path convex hull containing the origin must have smaller
dimension than its $j$-dimensional parameter simplex.
We will use the following particular consequence:
for $i+2$ consecutive path vectors of rank at least $i+1$,
the convex hull does not contain $0$.
This is what prevents an exception at the origin in the next
counting identity.

\subsection{Using one side of a subspace to count intersections}

For $0\le i<d$, define
\[
 H_i^+=\{u\in H_i:\ell_{i+1}(u)>0\}.
\]
This is the positive halfspace within $H_i$ specified by
$\ell_{i+1}$.
For a tuple $u_0,\ldots,u_i$ in $\mathcal V$, define
\begin{equation}\label{eq:alpha}
 \alpha_i(u_0,\ldots,u_i)=
 \begin{cases}
 1,&\conv\{u_0,\ldots,u_i\}\cap H_i^+\ne\varnothing,\\
 0,&\text{otherwise}.
 \end{cases}
\end{equation}
This is the only meaning of $\alpha_i$:
it records an intersection on one specified side of $H_i$.
No function $\alpha_d$ is needed.

\begin{lemma}\label{lem:intersection}
Let $0\le i<d$ and $u_0,\ldots,u_i\in\mathcal V$.
If these vectors are dependent, their convex hull has no nonzero
intersection with $H_i$.
If they are independent, any intersection with $H_i$ is a single
nonzero point, has strictly positive coefficients in the convex
combination, and satisfies $\ell_{i+1}(u)\ne0$.
Consequently,
\begin{equation}\label{eq:halves}
 \alpha_i(u_0,\ldots,u_i)+\alpha_i(-u_0,\ldots,-u_i)
 =
 \begin{cases}
 1,&\conv\{u_0,\ldots,u_i\}\cap(H_i\setminus\{0\})
       \ne\varnothing,\\
 0,&\text{otherwise}.
 \end{cases}
\end{equation}
\end{lemma}

\begin{proof}
If the vectors are dependent, their span has dimension at most
$i$.
By \eqref{eq:projection}, a vector in that span annihilated by
$\ell_1,\ldots,\ell_i$ must be zero.
Thus the span meets $H_i$ only at $0$.

Suppose they are independent.
Write a point of the convex hull as
$u=\sum_{s=0}^i t_su_s$, where $t_s\ge0$ and $\sum_st_s=1$.
For this point to lie in $H_i$, its coefficients must satisfy
\begin{equation}\label{eq:coefficients}
 \sum_st_s\ell_h(u_s)=0\quad(1\le h\le i),
 \qquad \sum_st_s=1.
\end{equation}
The homogeneous equations have a one-dimensional kernel by
\eqref{eq:projection}.
If a nonnegative solution exists, normalizing its coordinate
sum determines it uniquely.

None of its coefficients can vanish.
Otherwise $u$ would lie in the span of at most $i$ vectors.
Since $u\in H_i$, \eqref{eq:projection} would then force $u=0$.
This is impossible for a normalized combination of independent
vectors.
Finally, if $\ell_{i+1}(u)=0$, then $u\in H_{i+1}$, and
\eqref{eq:projection} would again force $u=0$.

Thus an intersection, when present, has exactly one of the two
signs for $\ell_{i+1}$.
Negating all vectors reverses that sign.
This proves \eqref{eq:halves}.
\end{proof}

For an $i$-face $F$, $0\le i<d$, let
\[
 f_i(F)=\sum_{P=(x_0,\ldots,x_i)\in\mathcal P(F)}
        \alpha_i\bigl(v(x_0),\ldots,v(x_i)\bigr).
\]
Write $\bar F=\{\bar x:x\in F\}$.
Taking antipodes and reversing the vertex order sends
$\mathcal P(F)$ bijectively to $\mathcal P(\bar F)$.
Since $v(\bar x)=-v(x)$ and the convex hull does not depend
on vertex order, \eqref{eq:halves} gives
\begin{equation}\label{eq:count-halves}
 c_i=\sum_{F\in X_i}\bigl(f_i(F)+f_i(\bar F)\bigr)
 \qquad(0\le i<d).
\end{equation}
This identity explains why antipodal symmetry is useful:
counting both positive sides is exactly counting every
nonzero intersection with $H_i$.

\subsection{The local identity: counting the endpoints of an interval}

The following lemma is the local reason why the parity in
\eqref{eq:goal} is preserved.
A hat denotes deletion of a position; equal vectors at different
positions are still treated separately.

\begin{lemma}\label{lem:deletion}
Let $0\le i<d$, and let
$u_s=e_{a_s}-e_{b_s}$, $0\le s\le i+1$, where $a_s\ne b_s$
and consecutive pairs $(a_s,b_s)$ have the same first or
second entry.
Put $\varepsilon=1$ if
\[
 \conv\{u_0,\ldots,u_{i+1}\}\cap(H_{i+1}\setminus\{0\})
 \ne\varnothing,
\]
and put $\varepsilon=0$ otherwise.
Then
\begin{equation}\label{eq:deletion}
 \sum_{s=0}^{i+1}
 \alpha_i(u_0,\ldots,\widehat{u_s},\ldots,u_{i+1})
 \equiv\varepsilon\pmod2.
\end{equation}
\end{lemma}

\begin{proof}
Put $U=\spanof\{u_0,\ldots,u_{i+1}\}$.
If $\dim U\le i$, both sides vanish by \eqref{eq:projection}.
Otherwise Lemma~\ref{lem:rank} gives
$0\notin\conv\{u_0,\ldots,u_{i+1}\}$.
Define
\[
 M=(\ell_h(u_s))_{\substack{1\le h\le i+1\\0\le s\le i+1}},
 \qquad b=(0,\ldots,0,1)^T,
 \qquad
 \mathcal J=\{t\in\R^{i+2}:t\ge0,\ Mt=b\}.
\]
By \eqref{eq:projection}, $M$ has rank $i+1$.
Thus $Mt=b$ defines an affine line.

Every feasible point has at least $i+1$ positive coordinates.
Indeed, otherwise \eqref{eq:projection} would force
$\sum_st_su_s=0$, contradicting
$\ell_{i+1}(\sum_st_su_s)=1$.
The same argument shows that, if exactly $i+1$ coordinates
are positive, the corresponding vectors are independent.

Consequently, every finite endpoint of $\mathcal J$ has
exactly one zero coordinate.
Conversely, at a feasible point with one zero coordinate,
the corresponding square submatrix of $M$ is invertible by
\eqref{eq:projection}.
The affine line therefore meets that coordinate hyperplane
at exactly one point, which is an endpoint.
Moreover, $\mathcal J$ cannot be a singleton:
with at most one zero coordinate, a small move along the line
is possible in at least one direction.

Normalizing an endpoint's coefficients to sum to one gives
a deletion tuple counted by $\alpha_i$.
Conversely, dividing the coefficients of a convex combination
counted by $\alpha_i$ by its positive $\ell_{i+1}$-value
gives a feasible point with one zero coordinate.
Hence the left side of \eqref{eq:deletion} counts precisely
the finite endpoints of $\mathcal J$.

Since the full convex hull excludes zero, we have
$\varepsilon=1$ if and only if $\ker M$ contains a nonzero
nonnegative vector $c$.
Every such $c$ is strictly positive.
Otherwise it uses at most $i+1$ vectors, so
\eqref{eq:projection} and $Mc=0$ imply
$\sum_sc_su_s=0$, contradicting the exclusion of zero
from the convex hull.

If $\varepsilon=1$, adding sufficiently large multiples
of $c$ to any solution of $Mt=b$ makes all coordinates
positive.
Since $\ker M$ is one-dimensional, $\mathcal J$ is therefore
a ray and has one finite endpoint.

If $\varepsilon=0$, every nonzero vector in $\ker M$
has both positive and negative coordinates.
Thus $\mathcal J$ is empty or bounded.
Being a closed convex subset of a line and never a singleton,
it is empty or a nondegenerate closed segment, with zero
or two endpoints.
The number of finite endpoints is therefore odd exactly
when $\varepsilon=1$, proving \eqref{eq:deletion}.
\end{proof}

\section{The contradiction}
\label{local}

A facet of a $j$-face is a contained $(j-1)$-face.

\begin{proposition}\label{prop:face}
For every $j$-face $F$, where $1\le j\le d$,
\begin{equation}\label{eq:face}
 \sum_{G\text{ facet of }F}f_{j-1}(G)
 \equiv
 \#\{P\in\mathcal P(F):
       S(P)\cap(H_j\setminus\{0\})\ne\varnothing\}
 \pmod2.
\end{equation}
\end{proposition}

\begin{proof}
Apply Lemma~\ref{lem:deletion}, with $i=j-1$, to each monotone
path in $F$ and sum.
Consecutive labels have a common entry, as required by the lemma.

On the left, deleting an internal vertex is paired with exchanging
the two coordinate changes immediately before and after that
vertex.
The paired paths have exactly the same vertex sequence
after deletion, so their contributions cancel modulo two.

Deleting the first vertex leaves a monotone path in a facet whose
new fixed coordinate is $1$.
Deleting the last leaves one in a facet whose new fixed coordinate
is $0$.
Every monotone path in a facet extends uniquely in the
corresponding way.
Thus precisely the left side of \eqref{eq:face} remains.
The right side is the number of nonzero intersections specified
in the proposition.
\end{proof}

It remains to check that the selected families $X_i$ have the
required boundary counts.

\begin{lemma}\label{lem:faces}
For $1\le j\le d$, a $(j-1)$-face $G$ is a facet of an odd number
of faces from $X_j$ if and only if $G\in X_{j-1}$ or
$\bar G\in X_{j-1}$.
\end{lemma}

\begin{proof}
Read the fixed bits of $G$ as a binary string of length
$K-j+1\ge2$.
A containing face in $X_j$ is obtained by freeing
one of these coordinates.
We therefore count the deletions that
leave an alternating string starting with $0$.

If the original string alternates starting with $0$, only deleting
the last bit works.
If it alternates starting with $1$, only deleting the first works.
Otherwise, a successful deletion must remove a member of the
unique equal adjacent pair; deleting either member leaves
the same string.
Thus in this case there are either zero or two successful
deletions.

The count is odd exactly for the two alternating strings, which
describe $G\in X_{j-1}$ and $\bar G\in X_{j-1}$.
\end{proof}

\begin{proof}[Proof of Theorem~\ref{thm:main}]
For $1\le j\le d$, sum Proposition~\ref{prop:face}
over $F\in X_j$.
By the definition of $c_j$,
\[
 c_j\equiv
 \sum_{F\in X_j}\sum_{G\text{ facet of }F}f_{j-1}(G)
 \pmod2.
\]
Lemma~\ref{lem:faces} identifies exactly the terms with odd
multiplicity.
Hence
\[
 c_j\equiv
 \sum_{G\in X_{j-1}}
 \bigl(f_{j-1}(G)+f_{j-1}(\bar G)\bigr)
 =c_{j-1}\pmod2,
\]
where the final equality is \eqref{eq:count-halves}.
This proves \eqref{eq:goal} and completes the proof
of Theorem~\ref{thm:main}.
\end{proof}

\end{document}